\documentclass[sn-mathphys-num]{sn-jnl}
\usepackage{graphicx}
\usepackage{multirow}
\usepackage{amsmath,amssymb,amsfonts}
\usepackage{amsthm}
\usepackage{mathrsfs}
\usepackage[title]{appendix}
\usepackage{xcolor}
\usepackage{textcomp}
\usepackage{manyfoot}
\usepackage{booktabs}
\usepackage{algorithm}
\usepackage{algorithmicx}
\usepackage{algpseudocode}
\usepackage{listings}
\newtheorem{example}{Example}
\theoremstyle{case}

\usepackage{algorithm}
\usepackage{algpseudocode}
\usepackage{graphicx}
\usepackage{subcaption}
\usepackage{caption}
\theoremstyle{thmstyleone}

\theoremstyle{thmstyletwo}

\theoremstyle{thmstylethree}

\begin{document}
\title[{\bf Development and Analysis of Chien-Physics-Informed Neural Networks for Singular Perturbation Problems}]{\bf Development and Analysis of Chien-Physics-Informed Neural Networks for Singular Perturbation Problems}

\author*[1]{\fnm{Gautam} \sur{Singh}}\email{gautam@nitt.edu}
	
\author[1]{\fnm{Sofia} \sur{Haider}}\email{Sofia.2025.nitt@gmail.com}
	
\affil[1]{\orgdiv{Department of Mathematics}, \orgname{National Institute of Technology Tiruchirappalli}, \city{Tiruchirappalli}, \postcode{620015}, \state{Tamil Nadu}, \country{India}}

\abstract{In this article, we employ Chien-Physics Informed Neural Networks (C-PINNs) to obtain solutions for singularly perturbed convection-diffusion equations, reaction-diffusion equations, and their coupled forms in both one and two-dimensional settings. While PINNs have emerged as a powerful tool for solving various types of differential equations, their application to singular perturbation problems (SPPs) presents significant challenges. These challenges arise because a small perturbation parameter multiplies the highest-order derivatives, leading to sharp gradient changes near the boundary layer. To overcome these difficulties, we apply C-PINNs, a modified version of the standard PINNs framework, which is specifically designed to address singular perturbation problems. Our study shows that C-PINNs provide a more accurate solution for SPPs, demonstrating better performance than conventional methods.}

\keywords{ Physics-informed neural network, Singularly perturbed problem, Petrov-Galerkin method, mean squared error, CPINN.}
	
\pacs[MSC Classification]{65L10, 65L11, 65L20, 65L60, 65L70, 68T07.}
	
\maketitle	
\section{Introduction}
Many real-world systems, from fluid dynamics to chemical reactions, are governed by singularly perturbed differential equations (SPDEs). These equations are challenging to solve because a small perturbation parameter multiplies the highest-order derivatives, creating sharp variations near boundary layers that traditional numerical methods struggle to capture \cite{CENZhongdi, Fujun, JAA}. Classical approaches such as the Finite Difference Method (FDM) and Finite Element Method (FEM) typically require extreme mesh refinement in these regions, which significantly increases computational cost \cite{PINNFEM}.

Deep learning-based approaches, particularly Physics-Informed Neural Networks (PINNs), have recently emerged as powerful tools for solving differential equations \cite{goodfellow2016, alAradi2019}. PINNs employ neural networks to approximate solutions while enforcing the governing equations through a residual-based loss function \cite{SenWang, Cpinn}, ensuring that the predicted solution is consistent with the underlying physics \cite{rathish1, rathish2}. The loss function incorporates the residuals of the PDE, boundary conditions, and initial conditions, with derivatives computed via automatic differentiation \cite{goodfellow2016, alAradi2019, Cpinn}. However, when applied to SPDEs, standard PINNs often fail to resolve sharp gradients near boundary layers, resulting in significant errors \cite{yadav, JAA}.

To address this limitation, researchers have introduced Chien-PINNs (C-PINNs) \cite{Cpinn}, a modification of the standard PINN framework specifically designed for singular perturbation problems (SPPs) \cite{JAA}. The central idea draws from composite asymptotic expansions, where the solution is decomposed into two parts: a smooth, slowly varying component and a rapidly changing boundary-layer component. Instead of relying on a single neural network to approximate the entire solution, C-PINNs employ multiple specialized subnetworks: some dedicated to smooth regions and others to boundary layers. These subnetworks are then blended using exponential weighting functions, yielding a uniformly valid approximation across the domain.

In this work, we extend the C-PINN framework obtain solutions for singularly perturbed convection-diffusion equations, reaction-diffusion equations, and their coupled forms in both one and two-dimensional settings. Such equations frequently arise in applications including chemical transport, biological diffusion, and fluid dynamics. By leveraging the C-PINN methodology, we aim to develop an accurate and computationally efficient approach for solving these challenging problems.

Solving SPDEs has long been a focus of numerical analysis. Zhongdi \cite{CENZhongdi} developed an FDM-based scheme for coupled convection-diffusion equations, while adaptive mesh refinement has been applied to reaction-diffusion systems \cite{Fujun}. Although effective, these approaches become computationally expensive as problem complexity and dimensionality increase \cite{PINNFEM}. More recently, deep learning methods-particularly PINNs-have gained attention as alternatives to classical numerical solvers \cite{alAradi2019}. While PINNs provide flexibility by embedding governing equations directly into neural networks, they struggle with singular perturbation problems, failing to accurately capture boundary layers \cite{yadav, JAA}.

To overcome these difficulties, enhanced PINN architectures have been developed. General-Kindered PINNs (GK-PINNs), introduced by Wang et al. \cite{SenWang}, adaptively adjust the network structure to improve accuracy. C-PINNs, proposed by Wang et al. \cite{Cpinn}, employ composite asymptotic expansions to handle sharp gradients more effectively. Opschoor et al. \cite{JAA} demonstrated that C-PINNs outperform both standard PINNs and FEMs in singular perturbation settings. Similarly, Kumar et al. \cite{rathish1, rathish2} applied neural networks to convection–diffusion problems in porous media, showcasing the potential of deep learning in modeling complex physical systems.

Despite these advancements, PINN-based approaches remain relatively unexplored for coupled convection-diffusion and reaction-diffusion systems. These systems are particularly challenging due to the presence of multiple interacting boundary layers. In this work, we extend the C-PINN framework to handle such coupled SPDEs. By combining composite asymptotic expansions with a multi-subnetwork architecture, our approach captures both smooth and sharp solution components, maintaining high accuracy even in the presence of multiple boundary layers.

The structure of the paper is as follows: Section \ref{Maths} introduces the mathematical formulation of the coupled singularly perturbed system. Section 2 describes the C-PINN framework, including the composite asymptotic expansion approach and implementation details. Section 3 presents numerical results for coupled convection–diffusion and reaction–diffusion problems. Finally, Section 4 concludes the paper and outlines future research directions.

\section{Mathematical Framework and Implementation}\label{Maths}
	
\subsection{General Form}\label{gen}

Consider a general singularly perturbed boundary layer problems defined in the simulation domain $\Omega \subset  \mathbb{R}^d$ with its boundary denoted as $\partial \Omega$, given by:
\begin{align}\label{eq1}
	\displaystyle & \mathcal{D}_\epsilon \Big( \mathbf{u}(x,t), x, t, f; \frac{\partial \mathbf{u}}{\partial t}, \frac{\partial^2 \mathbf{u}}{\partial t^2}, \dots, \frac{\partial \mathbf{u}(x,t)}{\partial x}, \frac{\partial^2 \mathbf{u}}{\partial x^2}, \dots \Big) = 0,  \quad \forall (x,t) \in \mathcal{U}, 
%	\\ 
%	\displaystyle & \mathcal{D}_\mu \Big( \mathbf{u}(x,t) , x, t, f_2; \frac{\partial \mathbf{u}(x,t)}{\partial t}, \frac{\partial^2 \mathbf{u}}{\partial t^2}, \dots, \frac{\partial \mathbf{u}(x,t)}{\partial x}, \frac{\partial^2 \mathbf{u}}{\partial x^2}, \dots \Big) = 0,  \quad \forall (x,t) \in \mathcal{U}.
\end{align}
with boundary conditions 
\begin{align} \label{eq2}
	\displaystyle    \mathcal{B} \Big( \mathbf{u}(x,t), x, t, g; \frac{\partial \mathbf{u}}{\partial x}, \dots \Big)  = 0, \; \quad \forall (x,t) \in \partial \mathcal{U},
\end{align}
and initial conditions 
\begin{align}\label{eq3} 
	\displaystyle    \mathcal{I}\Big( \mathbf{u}(x,t), x, t, h; \frac{\partial \mathbf{u}}{\partial t}, \dots \Big) = 0, \;  \quad \forall (x,t) \in \Gamma.
\end{align}
Where \( \mathbf{u} \) is a vector-valued function defined as:
$ \displaystyle \mathbf{u}(x,t): \mathbb{R}^{d+1} \to \mathbb{R}. $
The term $ \mathcal{D_{\epsilon}} $ denotes the residual form of the corresponding equations including source function $f(x,t)$, contains the differential operators with a small perturbation parameters $ \displaystyle \epsilon \in (0,1] $ multiplying the highest-order derivative. The term $ \mathcal{B} $ represent the boundary conditions in their residual form, including  source functions $ g(x,t)$. In the same way, $ \mathcal{I}$ represent the initial conditions in their residual form, including  source functions $ h(x,t).$ 
We define the space-time domain and its boundaries as follows:
$ \mathcal{U} = \{(x,t) \mid x \in \Omega, t \in [0,T]\},  \mathcal{\partial U} = \{(x,t) \mid x \in  \partial \Omega, t \in [0,T]\} $ \text{and}  $ 
\Gamma  =\{(x,t) \mid x \in  \partial \Omega, t =0\} . $

\subsection{Composite Asymptotic Expansion}
	SPPs have small parameters like $\epsilon$, which make the solution behave differently in different parts of the domain. There is a smooth solution (outer solution) and a rapidly changing solution near boundaries (inner solution). Traditional asymptotic methods can approximate these solutions separately, but they often fail to fully capture the behavior of the system. To address this, composite asymptotic expansion combines both the inner and outer solutions, creating a single, uniform approximation that works across the entire domain. 
	For a coupled system, we can write the solutions for \( u(x) \)  as:
	\begin{equation}\label{eq4}
		u^{\text{comp}}(x,t) =
		\sum_{k=0}^{n} \epsilon^k u_{k}^{1}(x,t) +
		\sum_{k=0}^{n} \epsilon^k u_{k}^{2}(x,t) 
		\left( e^{-p_1(x)/\delta_1} + e^{-p_2(x)/\delta_1} \right),
	\end{equation}
	
%	\begin{equation}
%		u_2^{\text{comp}}(x,t) =
%		\sum_{k=0}^{n} \mu^k u_{2,k}^{1}(x,t) +
%		\sum_{k=0}^{n} \mu^k u_{2,k}^{2}(x,t) 
%		\left( e^{-p_1(x)/\delta_2} + e^{-p_2(x)/\delta_2} \right).
%	\end{equation}
	
	Here, \( p_1(x) = x \) and \( p_2(x) = 1 - x \) represent the boundary layers at \( x = 0 \) and \( x = 1 \), while \( \delta_1 \) is the thicknesses of this layer.
	
\subsection{C-PINN Method}

C-PINN method is build on the foundation of composite asymptotic expansions, to solve singularly perturbed boundary-layer problems. The approximate solution $\mathbf{u}$ is constructed using multiple parallel sub-neural networks: outer NNs $\mathbf{u}$ to capture the smooth part in the solution and inner NNs $\mathbf{u}_{j}$ to capture the sharp gradients near boundary layers in the solution, here, $j$ =0,1. These sub-solutions are then integrated using an exponential weighting scheme, $ \exp\big(-\frac{p(x)}{\partial(\epsilon)}\big)$, ensuring a smooth transition between different solution regions. This entire NN framework is commonly known as the composite NN and represent a uniformly valid solution. The solution to equations (\ref{eq1})- (\ref{eq4}) is computed via optimization of the parameters of this composite NN using a suitable optimizer.
The composite neural network with input $(x,t)$ and output $u$ is given as: 
	\begin{align}\label{eq5}
			u(x,t; \theta_{1}, \theta_{2},\theta_{3}) = u(x,t;\theta_{1}) + u_{0}(x,t;\theta_{2})*e^{-p_1(x)/\delta_1} + u_{1}(x,t;\theta_{3})*e^{-p_2(x)/\delta_1} 
	\end{align}
	
where $u $ is the outer NNs and $u_{0}$, $u_{1}$ are inner NNs and $\theta_{i}, i $ = 1 to 3 denote the parameters of these networks.

\subsection{Algorithm of the C-PINN Method}
Now, we discuss the algorithm of the C-PINN method which is as follows:
Firstly, identify $p(x)$  and $\partial(\epsilon)$  to define the location and thickness of all the boundary layers present in equation. For example, for the boundary layer with a thickness $O (\epsilon^n   ) $  near $x = x_{0}$, we take $p(x) = ||x - x_{0}||_{2}$ and $\partial(\epsilon) = \epsilon^n$. 
For instance, in the case of a boundary layer with thickness of the order of $O(\epsilon^n)$ near the point $x = x_0$, we define the function $p(x)$ as the Euclidean norm $||x - x_0||_2$, and specify the layer width function as $\partial(\epsilon) = \epsilon^n$.
Secondly, we construct the composite NN as defined in equation (\ref{eq5}) based on the given problem. For example, in a convection-diffusion system where the boundary layer forms only at one end of the coupled system say, at $(x = 0)$, the network will include only $ u $ and $ u_{0} $. However, in a reaction-diffusion system, where boundary layers may appear at both the ends, the full structure described in equation (\ref{eq5}) represents the composite NN.
Lastly, the loss function is computed using automatic differentiation implemented in PyTorch, and is defined as follows:  
	\begin{align}
		L(\theta)  &= \lambda_{\mathcal{D_{\epsilon}}} L_{D_{\epsilon}}(\theta) + \lambda_{{B}} L_{B}(\theta)  + \lambda_{{I}} L_{I}(\theta)  .
	\end{align}
	where; $ \theta \equiv (\theta_1, \theta_2, \theta_3)$ and the coefficients  $\lambda_{\mathcal{D_{\epsilon}}}, \lambda_{{B}}, \lambda_{{I}}$ determine the contribution of individual loss term in the total loss. All the loss terms are defined as follows:
	\begin{align*}
		L_{D_{\varepsilon}}(\theta) &= \frac{1}{N_D} \sum_{i=1}^{N_D} \left\| \mathcal{D}_{\varepsilon} (u(x^{(i)}, t^{(i)}; \theta), x^{(i)}, t^{(i)}, f_1^{(i)}) \right\|_2^2, \\
%		L_{D_{\mu}}(\theta) &= \frac{1}{N_D} \sum_{i=1}^{N_D} \left\| \mathcal{D}_{\mu} (u(x^{(i)}, t^{(i)}; \theta), x^{(i)}, t^{(i)}, f_2^{(i)}) \right\|_2^2, \\
		L_B(\theta) &= \frac{1}{N_B} \sum_{i=1}^{N_B} \left\| \mathcal{B} (u(x^{(i)}, t^{(i)}; \theta), x^{(i)}, t^{(i)}, g^{(i)}) \right\|_2^2, \\
		L_I(\theta) &= \frac{1}{N_I} \sum_{i=1}^{N_I} \left\| \mathcal{I} (u(x^{(i)}, t^{(i)}; \theta), x^{(i)}, t^{(i)}, h^{(i)}) \right\|_2^2.
	\end{align*}
	where, the loss terms $ L_{D_{\varepsilon}}(\theta)$ represents the Mean Squared Error (MSE) of equations (\ref{eq1})-(\ref{eq2}), computed at the collocation points. Additionally, the loss terms $ L_B(\theta) $ and $ L_I(\theta) $ correspond to the MSE of the boundary and initial conditions in equations (\ref{eq3})-(\ref{eq4}), respectively. And thus, we conclude the process by minimizing the loss function \( \mathcal{L}(\theta) \) through the optimization of parameters  $ \theta \equiv (\theta_1, \theta_2, \theta_3)$ using a suitable optimizer such as L-BFGS, RMSProp, or Adam.
	
	\section{Computational Results and Analysis}\label{result}
	In this section, we explore the application of C-PINN to solve various types of SPPs. We begin with the convection-diffusion problem and reaction-diffusion problem followed by their coupled system and 2 Dimensional SPPs.
	Each sub-neural network in the C-PINN architecture consists of fully connected layers, with a hyperbolic tangent activation function applied at each layer. 
	 To train the model, we rely on the Adam optimizer with a suitably chosen learning rate. The model's accuracy is evaluated using the L$_2$  loss, calculated as the MSE. 
	We use the Xavier initialization technique to set up the network's weights and biases. For selecting training points, we  utilize the Latin Hypercube Sampling (LHS) method to ensure a well-distributed dataset.

	\subsection{SPPs in One Dimension}
	In this section we use the NN architecture having 3 fully connected layers with each layer having 50 and 100 neurons respectively. Also, we train the NN for 10,000 epochs and 600 uniformly distributed collocation points and utilize the Adam optimizer with a learning rate of $5*10^{-4}$ in order to minimize the loss.
\begin{algorithm}
\caption{C-PINN Algorithm for 1D Singularly Perturbed Problem}
\begin{algorithmic}[1]
			\State \textbf{Initialize:}
			\State \hspace{1em} Initialize collocation points $x \in [0,1]$; learning rate $\eta$; maximum iterations $M$
			\State \hspace{1em} Initialize neural network weights: $\theta^o$, $\{\theta^k\}_{k=1}^K$ for the solution $u$
			\State \hspace{1em} Identify boundary layer location and thickness to obtain $p(x)$, $\delta(\epsilon)$
			
			\State Fix perturbation parameter $\varepsilon = \varepsilon^j$
			\State Define \textbf{OuterNN} for capturing smooth solution
			\State Define \textbf{InnerNN}$_{0}$, \textbf{InnerNN}$_{1}$ for layers at $x=0$ and $x=1$
			
			\State Construct composite network $u^\theta(x)$:
			\Statex \hspace{1.5em} $u^\theta(x) \gets \text{OuterNN}(x) + \text{InnerNN}_{0}(x)e^{-x/\epsilon} + \text{InnerNN}_{1}(x)e^{-(1-x)/\epsilon}$
			
			\State Define residual operator:
			\Statex \hspace{1.5em} $R(x) \gets \mathcal{D}^\varepsilon[u^\theta](x) - f(x)$
			
			\For{$m = 1$ to $M$}
			\State Sample training points $x \in [0,1]$
			\State Compute $u^\theta(x)$ and its derivatives via auto-diff
			\State Evaluate residual $R(x)$ at collocation points
			\State Compute loss:
			\Statex \hspace{1.5em} $\mathcal{J} \gets \mathbb{E}[R(x)^2] + \text{MSE of boundary conditions}$
			\State Backpropagate gradients and update parameters:
			\Statex \hspace{1.5em} $\theta \gets \theta - \eta \nabla_\theta \mathcal{J}$
			\EndFor
			
			\State \textbf{Return:} Trained network $\Rightarrow$ composite solution $u^\theta(x)$
			
			\State Use current weights as initialization for next continuation step and repeat for $M$ iterations.
\end{algorithmic}
\end{algorithm}

	\subsubsection{Singularly Perturbed Convection-Diffusion Problem}
	Consider the following convection-diffusion problem:   
$$ \displaystyle  -\epsilon y''(x) + y'(x) + y(x) = f(x), \quad x \in [0,1],$$
where \( \epsilon \) is a small perturbation parameter and 
the boundary conditions are $ y(0) = 0\; \text{and} \; y(1) = 0. $ \\
The source term $f(x)$ is chosen such that the analytical solution  is given by: $$\displaystyle  y(x) = (1 - e^{\frac{x-1}{\epsilon}}) \sin(x). $$  
A NN \( y_{\theta}(x) \) is used to approximate \( y(x) \). To enforce the differential equation, we compute the residual:  
$\displaystyle  R(x) = -\epsilon y_{\theta}''(x) + y_{\theta}'(x) + y_{\theta}(x) - f(x). $ To improve stability, a soft residual-based weighting function is applied $\displaystyle  w(x) = e^{-\lambda |R(x)|}, $ where \( \lambda \) is a scaling parameter to control the weight decay. \\
The total loss function to minimize is: 
	\begin{align}
		\displaystyle 	\mathcal{L} &= \mathcal{L}_{\text{ODE}} + \mathcal{L}_{\text{BC}} \nonumber \\
		&= \sum_{x \in \Omega} w(x) R(x)^2 + (y_{\theta}(0) - 0)^2 + \lambda_{\text{BC}} (y_{\theta}(1) - 0)^2. \nonumber
	\end{align}
	The NN used in PINN approach consists of three fully connected layers, each with 150 neurons, and employs the \textit{tanh} activation function. Training is performed using the Adam optimizer with a learning rate of $ 10^{-3}, $ along with a step decay scheduler to gradually reduce the learning rate over time. The perturbation parameter is set to  $ \epsilon = 0.01$. To better enforce the boundary conditions, a weighting factor of  
	$ \lambda_{\text{BC}} = 3.0 $ is applied specifically near \( x = 1 \), where sharp gradients occur due to the singular perturbation. Additionally, a soft residual weighting factor of  $ \lambda = 0.8 $ is used to stabilize training by reducing the influence of large residuals.\\
	The model is trained for 30,000 epochs, achieving a low final loss, demonstrating that the NN effectively captures the solution to the problem. \\
	\begin{figure}[h]
		\centering
		\includegraphics [width=0.45\textwidth]{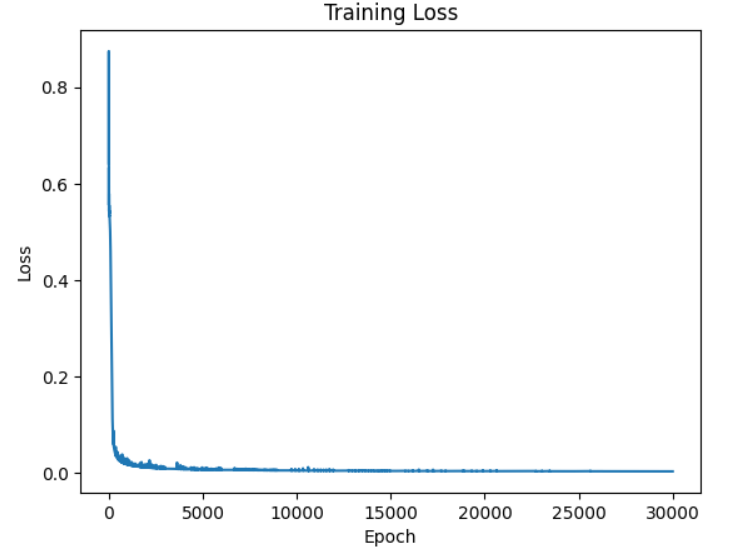}
		\includegraphics [width=0.45\textwidth]{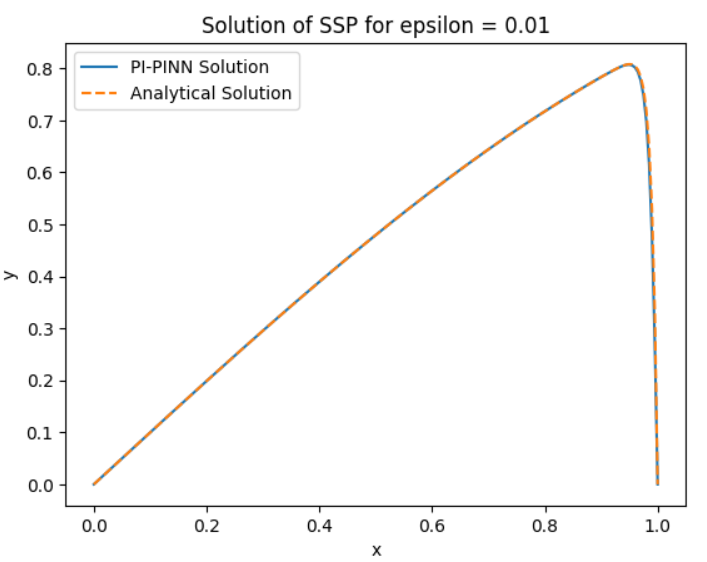}
		\caption{ Loss and solution plot for PI-PINN for different epochs.}
		\label{figure1}
	\end{figure}

	\begin{figure}[h]
		\centering
		\includegraphics [width=0.45\textwidth]{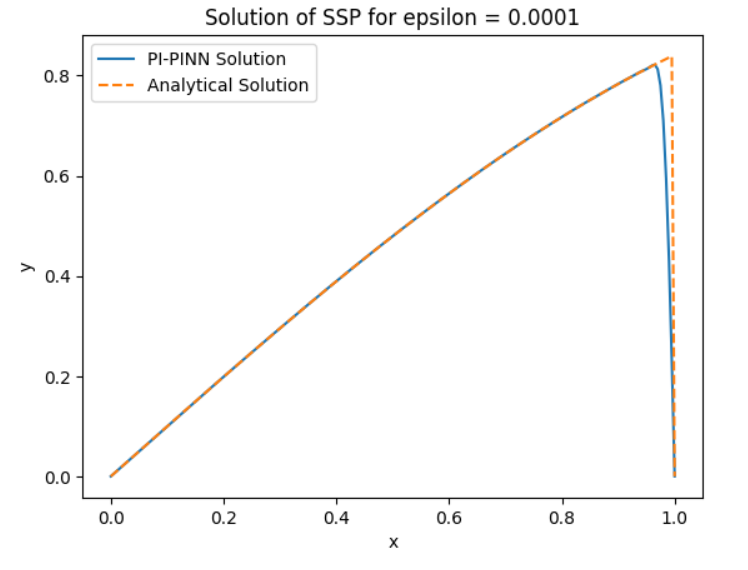}
		\includegraphics [width=0.45\textwidth]{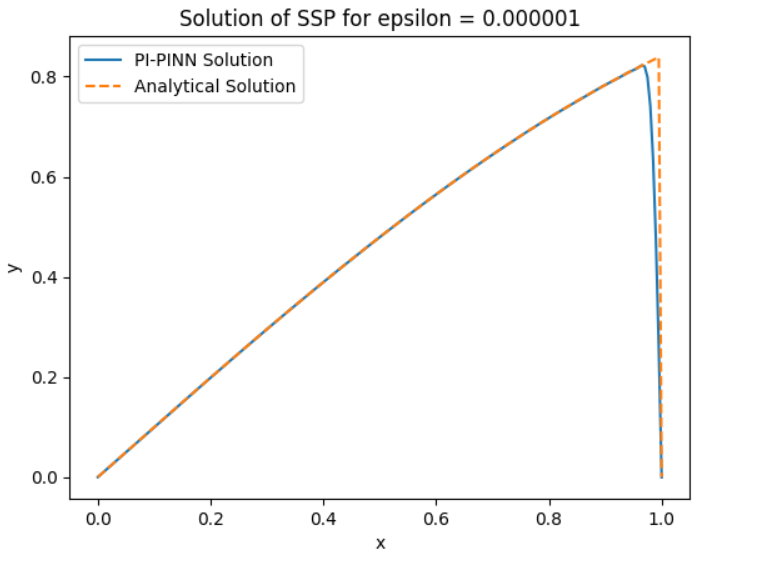}
		\caption{ Solution plot for PI-PINN for $\epsilon$= 0.0001 and $\epsilon$=0.000001.}
		\label{figure2}
	\end{figure}	
	
	But still this approach was not able to fully approximate the solution as $\epsilon$ decreased. So, we move forward with this problem and apply C-PINN approach and then compare with the PI-PINN approach which we used previously.
	In this problem, the boundary layer is located at $x=1$, with a thickness of order $\mathcal{O}(\epsilon).$ For our setup, we take $\epsilon = 10^{-5}$.
	Now, let's go through the step-by-step implementation of the C-PINN method and the training process of the model. \\
	First, we define two types of sub-neural networks: outer NNs and inner NNs. The outer networks capture the smooth variations, while the inner networks handle the sharp gradients of the problem. Each network consists of three fully connected layers with 50 neurons per layer. 
	We define the PI-PINN model and create model instances, namely PI-PINN, outer-nn, and inner-nn. Additionally, we use the following composite neural network to represent a uniform solution:
	{\footnotesize	\begin{equation}
			\displaystyle y(x;\theta_1,\theta_2)  =
			\text{outer\_nn}(x,\theta_1) + \text{inner\_nn}(x,\theta_2) \cdot \text{exp}\left(\frac{-(1-x)}{\epsilon}\right)
	\end{equation}}

	Next, we define the loss function for both C-PINN and PI-PINN using PyTorch's built-in automatic differentiation. The loss formulation also includes the residual, ensuring the network learns to satisfy the governing equations. \\
	We take 600 uniformly distributed collocation points for training and train both the models namely PI-PINN and C-PINN for 10000 epochs using Adam optimizer with learning rate $10^{-3}$ and $5*10^{-4}$ respectively.\\

	\begin{figure}[h]
		\centering
		\includegraphics [width=0.45\textwidth]{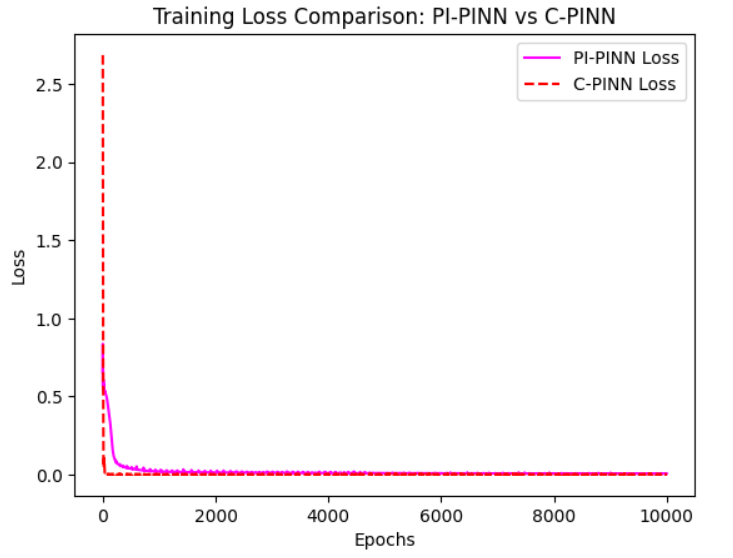}
		\includegraphics [width=0.45\textwidth]{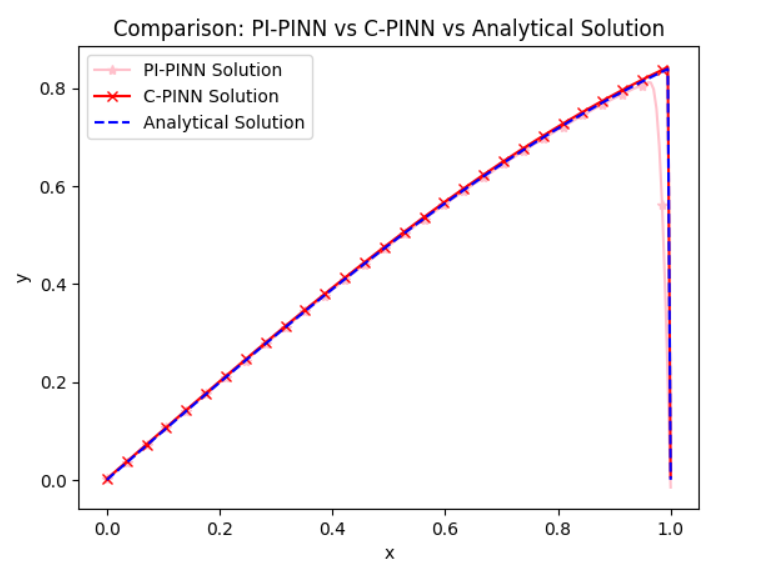}
		\caption{ Loss and solution plot for PI-PINN and C-PINN for different epochs.}
		\label{figure3}
	\end{figure}

	Table \ref{table1} shows that initially, C-PINN loss starts with a high loss of 2.686301, much greater than PI-PINN's 0.654495. However, C-PINN rapidly decreases to 0.000123 by epoch 500 and further stabilizes around 0.000013 at epoch 9500. While, PI-PINN starts lower but decreases more gradually, from 0.654495 to 0.006071 over the same period. Moreover, consistently achieves a lower loss compared to PI-PINN across all epochs, indicating better convergence and accuracy. This comparison highlights the efficiency and stability of C-PINN in solving the given problem.
	
	\begin{table}[ht]
		\centering
		\renewcommand{\arraystretch}{1.4} % Adjust row height
		\setlength{\tabcolsep}{4pt} % Adjust column spacing
		\tiny % Reduce font size for better fit
		\begin{tabular}{|c|c | c |c |c |c |c |c |c |c| c|}
			\hline
			\textbf{Epoch} & 0 & 500 & 1000 & 1500 & 2000 & 2500 & 3000 & 3500 & 4000 & 4500  \\ 
			\hline
			\textbf{PI-PINN Loss} & 0.654495 & 0.037589 & 0.018725 & 0.015025 & 0.013057 & 0.011936 & 0.011880 & 0.009979 & 0.009359 & 0.008878 \\ 
			\textbf{C-PINN Loss}  & 2.686301 & 0.000123 & 0.000054 & 0.000043 & 0.000041 & 0.000034 & 0.000033 & 0.000029 & 0.000026 & 0.001955 \\ 
			\hline
			\textbf{Epoch} & 5000 & 5500 & 6000 & 6500 & 7000 & 7500 & 8000 & 8500 & 9000 & 9500  \\ 
			\hline
			\textbf{PI-PINN Loss} & 0.008308 & 0.009860 & 0.009854 & 0.007913 & 0.007227 & 0.006645 & 0.006418 & 0.006445 & 0.010089 & 0.006071 \\ 
			\textbf{C-PINN Loss}  & 0.000019 & 0.000026 & 0.000015 & 0.000050 & 0.000032 & 0.000016 & 0.000012 & 0.000009 & 0.000144 & 0.000013 \\ 
			\hline
		\end{tabular}
		\caption{\( L_2 \) loss for PI-PINN and C-PINN at different epochs.}
		\label{table1}
	\end{table}
	
	\subsubsection{Singularly Perturbed Reaction-Diffusion Problem}
	Consider the following reaction-diffusion problem:  
	$$ \displaystyle  -(\epsilon)^2 u''(x) + 8 u(x) = f(x), \quad x \in [0,1]	$$
	with boundary conditions $ u(0) = 0$\; \text{and} $u(1) = 0$, where \( \epsilon \) is a small perturbation parameter.  \\
	The source term is chosen such that the analytical solution of the equation is given by:  
	$$\displaystyle  u(x) = e^{\frac{-x}{\epsilon}} +  e^{\frac{-(1-x)}{\epsilon}} -1 - e^{-\frac{1}{\epsilon}}. $$  
	In this problem, the boundary layer is located at both the ends at $x=0$ and $x=1$, with a thickness of order $\mathcal{O}(\epsilon).$ For our setup, we take $\epsilon = 10^{-5}$.
	After we define two types of sub-neural networks: outer NNs and inner NNs, we utilize the following composite NN to represent a uniform solution: 
    {\scriptsize   
    	\begin{align}
    		\displaystyle u(x;\theta_1,\theta_2, \theta_3) =
    		\text{outer\_nn}(x,\theta_1) + \text{inner1\_nn}(x,\theta_2) \cdot \text{exp}\left(\frac{-x}{\epsilon}\right) + \text{inner2\_nn}(x,\theta_3) \cdot \text{exp}\left(\frac{-(1-x)}{\epsilon}\right)
    	\end{align}
    }

	Table \ref{table2} shows that initially, C-PINN starts with a high loss of 28.686800. However, C-PINN rapidly decreases to 0.038382 by epoch 500 and further stabilizes around 0.003206 at epoch 9500. C-PINN reaches a loss of order  $10^{-3}$ in the final loss. C-PINN's sharp loss decline shows it captures the solution more accurately and adapts better, making it an effective approach.\\
	\begin{figure}[h]
		\centering
		\includegraphics [width=0.45\textwidth]{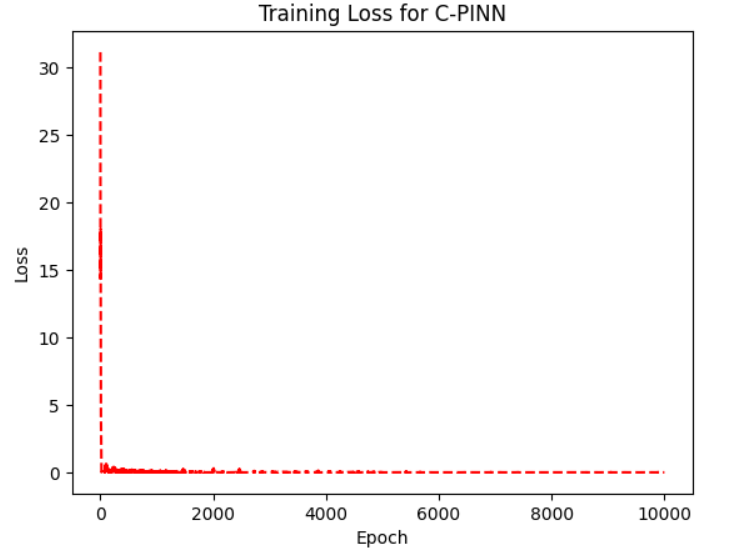}
		\includegraphics [width=0.45\textwidth]{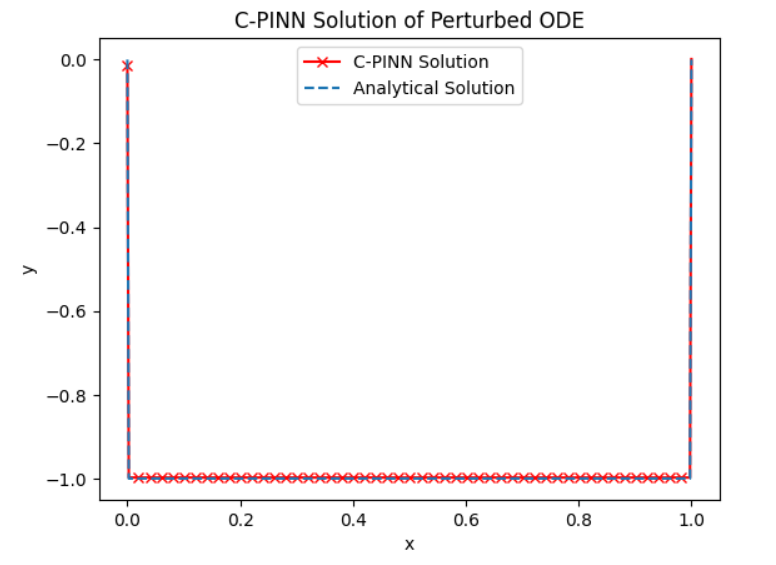}
		\caption{ Loss and solution plot for PI-PINN and C-PINN for different epochs.}
		\label{figure4}
	\end{figure}	
	
	\begin{table}[ht]
		\small
		\centering
		\renewcommand{\arraystretch}{1.4} % Adjust row height
		\setlength{\tabcolsep}{4pt} % Adjust column spacing
		\tiny % Reduce font size for better fit
		\begin{tabular}{|c|c | c |c |c |c |c |c |c |c| c|}
			\hline
			\textbf{Epoch}  & 0 & 500 & 1000 & 1500 & 2000 & 2500 & 3000 & 3500 & 4000 & 4500 \\ 
			\hline
			\textbf{C-PINN Loss} & 28.686800 & 0.038382 & 0.069121 & 0.021204 & 0.003783 & 0.014811 & 0.012587 & 0.011475 & 0.020970 & 0.004023 \\ 
			\hline
			\textbf{Epoch}  & 5000 & 5500 & 6000 & 6500 & 7000 & 7500 & 8000 & 8500 & 9000 & 9500 \\ 
			\hline
			\textbf{C-PINN Loss} & 0.007388 & 0.003242 & 0.003943 & 0.006221 & 0.003254 & 0.004066 & 0.003756 & 0.004996 & 0.003793 & 0.003206 \\ 
			\hline
		\end{tabular}
		\caption{\( L_2 \) loss for C-PINN at different epochs.}
		\label{table2}
	\end{table}

\begin{algorithm}
	\caption{C-PINN Algorithm for Coupled System}
	\begin{algorithmic}[1]
		\State \textbf{Initialize:}
		\State \hspace{1em} Initialize collocation points $x \in [0,1]$; learning rate $\eta$; maximum iterations $M$
		\State \hspace{1em} Initialize neural network weights: $\theta_1^o$, $\{\theta_1^k\}_{k=1}^K$, $\theta_2^o$, $\{\theta_2^k\}_{k=1}^K$ for $u_1$ and $u_2$
		\State Figure out the boundary layer location and boundary layer thickness to find $p(x)$,  $\delta(\epsilon)$ 
	
		\State Fix $\varepsilon = \varepsilon^j$, $\mu = \mu^j$
				\State Define \textbf{OuterNN}$_1$, \textbf{OuterNN}$_2$ for smooth solution components
		\State Define \textbf{InnerNN}$_{10}$, \textbf{InnerNN}$_{11}$, \textbf{InnerNN}$_{20}$, \textbf{InnerNN}$_{21}$ for boundary layers
	\State Construct composite networks $u_1^\theta$, $u_2^\theta$ using outer and inner subnetworks
		\Statex \hspace{1.5em} $u_1^\theta(x) \gets \text{OuterNN}_1(x) + \text{InnerNN}_{10}(x)e^{-x/\epsilon} + \text{InnerNN}_{11}(x)e^{-(1-x)/\epsilon}$
		\Statex \hspace{1.5em} $u_2^\theta(x) \gets \text{OuterNN}_2(x) + \text{InnerNN}_{20}(x)e^{-x/\epsilon} + \text{InnerNN}_{21}(x)e^{-(1-x)/\epsilon}$
		\State Evaluate residuals at collocation points
		\State Define loss $\mathcal{L}$ using PDE residuals, boundary conditions
		\For{$m = 1$ to $M$}
		\State Sample training points $x$ in domain $[0,1]$
		\State Compute $u_1(x), u_2(x)$ and their derivatives using auto-diff
		\State Define residuals:
	\Statex \hspace{1.5em} $R_1(x) \gets \mathcal{D}_1^{\varepsilon}[u_1^\theta, u_2^\theta] - f_1(x)$
	\Statex \hspace{1.5em} $R_2(x) \gets \mathcal{D}_2^{\mu}[u_1^\theta, u_2^\theta] - f_2(x)$
	
		\State Compute loss: 
		\Statex \hspace{1.5em} $\mathcal{L} \gets \mathbb{E}[R_1(x)^2 + R_2(x)^2] + \text{MSE of boundary conditions}$
		\State Backpropagate gradients and update weights using optimizer to minimize the loss
				\State Perform gradient descent step: $\theta \gets \theta - \eta \nabla_\theta \mathcal{L}$
		\EndFor
		\State \textbf{Return:} Trained networks $\Rightarrow$ composite solutions $u_1^\theta(x), u_2^\theta(x)$

		\State Use current weights as initialization for next step and continue the process for $M$ iterations.

	\end{algorithmic}
\end{algorithm}

\subsection{Coupled Systems of SPPs}
We define two types of NNs namely Outer NN and Inner NN. The NN consists of three fully connected layers having 100 and 150 neurons per layer respectively for Outer NN and Inner NN. Also, we define a `safe' exponential function that clamps input values before applying the exponential operation, preventing overflow and underflow issue. Loss is minimization using the `Adam' optimizer with a learning rate of $5*10^{-4}$. We take 600 uniformly distributed collocation points and train our NN model  for 7000 epochs.
\subsubsection{Coupled Systems of Convection-Diffusion Problems}
The general system of two coupled singularly perturbed convection-diffusion problem is given by:
\begin{align}
	-\varepsilon u_1''(x) - a_1(x) u_1'(x) + b_{11}(x) u_1(x) + b_{12}(x) u_2(x) &= f_1(x), \\  
	-\mu u_2''(x) - a_2(x) u_2'(x) + b_{21}(x) u_1(x) + b_{22}(x) u_2(x) &= f_2(x),  
\end{align}
where $ \mathbf{u} =  (u_1 ,u_2 )^{T} ,\displaystyle  x \in (0,1) $ and $ a_{i}, i=1,2$ denotes convection term and  $\displaystyle b_{ij},  i,j=1,2$ represents the coupling term with the boundary conditions:
\begin{equation}
	u_1(0) = u_2(0) = u_1(1) = u_2(1) = 0.
\end{equation}
Here, the boundary layer forms on either side of the boundary. Consider  the following system of two coupled singularly perturbed convection-diffusion problems \cite{CENZhongdi}:
\begin{align}
	-\varepsilon u_1''(x) - u_1'(x) + 2u_1(x) - u_2(x) &= f_1(x), \\
	-\mu u_2''(x) - 2u_2'(x) + 4u_2(x) - u_1(x) &= f_2(x), 
\end{align}
subject to the boundary conditions:
$$u_1(0) = u_1(1) = u_2(0) = u_2(1) = 0. $$ 
where $\epsilon $ and $\mu$ are  small perturbation parameters. Without loss of generality (WLOG), we assume $0<\mu < \epsilon \ll 1$.\\
The source terms  $ f_1(x), f_2(x) $ are chosen such that the system has the following analytical solution:
\begin{align*}
	u_1(x) &= \frac{1 - \exp(-x/\varepsilon)}{1 - \exp(-1/\varepsilon)} + \frac{1 - \exp(-x/\mu)}{1 - \exp(-1/\mu)} - 2 \sin \frac{\pi}{2} x, \\
	u_2(x) &= \frac{1 - \exp(-x/\mu)}{1 - \exp(-1/\mu)} - x \exp(x - 1).
\end{align*}
In this coupled system, the boundary layer is formed at $x=0$, with a thickness of order $\mathcal{O}(\epsilon).$ For our setup, we take $\epsilon = 10^{-7}$ and $\mu = 10^{-5}$.
We create model instances, namely outer-nn-1, outer-nn-2, inner-nn-10 and inner-nn-20 using the Outer NN and Inner NN. Moreover, we utilize the following composite NN to represent a uniform solution: 
\begin{align}\nonumber
	\displaystyle u(x,\theta_1,\theta_2,\theta_3,\theta_4) & =
	\text{outer\_nn}_1(x,\theta_1) + \text{inner\_nn}_{10}(x,\theta_2) \cdot \text{safe\_exp}\left(\frac{-x}{\epsilon}\right)  \\  \nonumber
	& +  \text{outer\_nn}_2(x,\theta_3) + \text{inner\_nn}_{20}(x,\theta_4) \cdot \text{safe\_exp}\left(\frac{-x}{\mu}\right) 
\end{align}

\begin{figure}[h]
	\centering
	\includegraphics [width=0.45\textwidth]{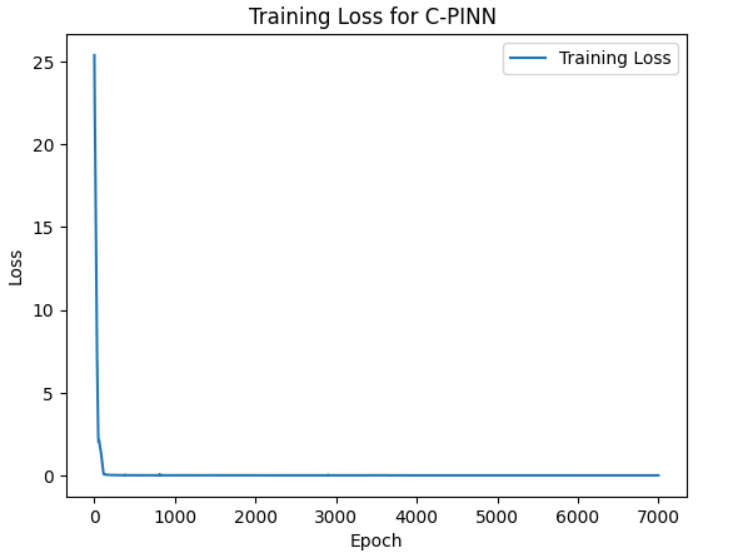}
	\includegraphics [width=0.45\textwidth]{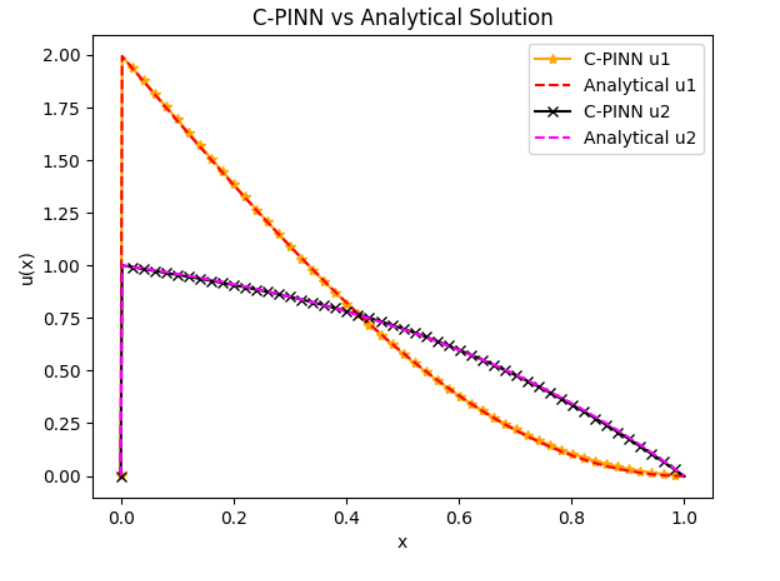}
	\caption{ Loss and solution plot for C-PINN for different epochs.}
	\label{figure5}
\end{figure}

\begin{figure}[h]
	\centering
	\includegraphics [width=0.47\textwidth]{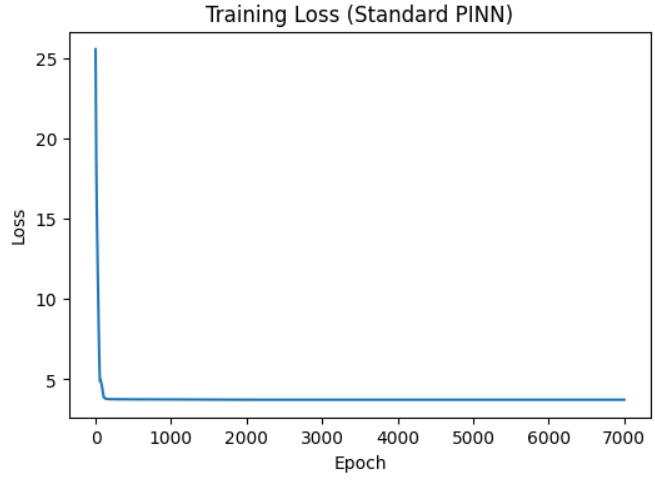}
	\includegraphics [width=0.44\textwidth]{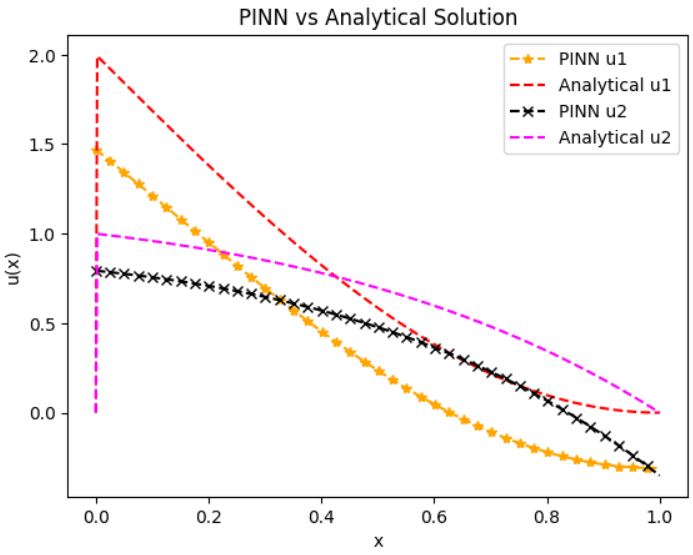}
	\caption{Loss and solution plot for PINN for different epochs.}
	\label{figure6}
\end{figure}	

\begin{table}[ht]
	\footnotesize
	\centering
	\renewcommand{\arraystretch}{1}
	\begin{tabular}{|c|c |c| c| c| c| c| c|}
		\hline
		\textbf{Epoch} & 0 & 500 & 1000 & 1500 & 2000 & 2500 & 3000 \\
		\hline
		\textbf{C-PINN Loss} 
		& 25.850483 & 0.015614 & 0.014740 & 0.013737 & 0.012002 & 0.011790 & 0.010592 \\
		\textbf{PINN Loss} 
		& 25.570000 & 3.760100 & 3.753100 & 3.745900 & 3.734800 & 3.734100 & 3.733900 \\
		\hline
		\textbf{Epoch} & 3500 & 4000 & 4500 & 5000 & 5500 & 6000 & 6500 \\
		\hline
		\textbf{C-PINN Loss} 
		& 0.009760 & 0.009230 & 0.008909 & 0.008491 & 0.008069 & 0.007949 & 0.007749 \\
		\textbf{PINN Loss} 
		& 3.733800 & 3.733800 & 3.733700 & 3.733700 & 3.733700 & 3.733600 & 3.733600 \\
		\hline
	\end{tabular}
	\caption{Comparison of C-PINN loss and PINN loss at different epochs.}
	\label{table4}
\end{table}

Table \ref{table4} shows that initially, both PINN and C-PINN start with comparable loss. However, C-PINN loss rapidly decreases to 0.015614 by epoch 500 and further stabilizes around 0.007749 at epoch 6500. Moreover, consistently achieves a lower loss compared to PINN across all epochs, indicating better convergence and accuracy. This comparison highlights the efficiency and stability of C-PINN in solving the given problem.\\

\subsubsection{Coupled Systems of Reaction-Diffusion Problems}
The general system of two coupled singularly perturbed reaction-diffusion problem is given by:
\begin{align}
	-\varepsilon^2 u_1''(x) + a_1(x) u_1(x) + b_{11}(x) u_1(x) + b_{12}(x) u_2(x) &= f_1(x), \\  
	-\mu^2 u_2''(x) + a_2(x) u_2(x) + b_{21}(x) u_1(x) + b_{22}(x) u_2(x) &= f_2(x),  
\end{align}
where $ \mathbf{u} =  (u_1 ,u_2 )^{T}  ,\displaystyle  x \in (0,1) $ and $ a_{i}$ where $i=1,2$ denotes reaction term and  $\displaystyle b_{ij}$ where  $i,j=1,2$ represents the coupling term with the boundary conditions:
\begin{equation}
	u_1(0) = u_2(0) = u_1(1) = u_2(1) = 0.
\end{equation}
Here, the boundary layers may appear on both sides of the boundary. \\
Consider  the following system of two coupled singularly perturbed reaction-diffusion problems \cite{gautam1}:
\begin{align}
	-\varepsilon^2 u_1''(x)  + 2u_1(x) - u_2(x) &= f_1(x), \\
	-\mu^2 u_2''(x)  - u_1(x) + 4u_2(x)  &= f_2(x), 
\end{align}
subject to the boundary conditions:
$$u_1(0) = u_1(1) = u_2(0) = u_2(1) = 0. $$ 
where $\epsilon $ and $\mu$ are  small perturbation parameters. WLOG, we assume $0<\mu < \epsilon \ll 1$.\\
The source terms  $ f_1(x), f_2(x) $ are chosen such that the system has the following analytical solution:
\begin{align*}
	u_1(x) &= \frac{ \exp(-x/\varepsilon) + \exp(-(1-x)/\varepsilon) }{1 - \exp(-1/\varepsilon)} + \frac{\exp(-x/\mu) + \exp(-(1-x)/\mu) }{1 - \exp(-1/\mu)} - 2 , \\
	u_2(x) &= \frac{ \exp(-x/\mu) + \exp(-(1-x)/\mu) }{1 - \exp(-1/\mu)} - 1 .
\end{align*}
In this coupled system, the boundary layer is present at both ends, specifically at $x=0$ and $x=1$, with a thickness of order $\mathcal{O}(\epsilon).$ For our setup, we take $\epsilon = 10^{-10}$ and $\mu = 10^{-8}$.
We create model instances, namely outer-nn-1, outer-nn-2, inner-nn-10, inner-nn-11, inner-nn-20 and inner-nn-21 and utilize the following composite NN to represent a uniform solution:

{\scriptsize \begin{align}\nonumber
		\displaystyle u(x,\theta_1,\theta_2,\theta_3,\theta_4, \theta_5, \theta_6) & =
		\text{outer\_nn}_1(x,\theta_1) + \text{inner\_nn}_{10}(x,\theta_2) \cdot \text{safe\_exp}\left(\frac{-x}{\epsilon}\right)+  \\ \nonumber
		&   \text{inner\_nn}_{11}(x,\theta_3) \cdot \text{safe\_exp}\left(\frac{-(1-x)}{\epsilon}\right)  
		+  \text{outer\_nn}_2(x,\theta_4)+ \\  \nonumber
		& \text{inner\_nn}_{20}(x,\theta_5) \cdot \text{safe\_exp}\left(\frac{-x}{\mu}\right)  
		+ \text{inner\_nn}_{21}(x,\theta_6) \cdot \text{safe\_exp}\left(\frac{-(1-x)}{\mu}\right). 
\end{align}}

\begin{figure}[h]
	\centering
	\includegraphics [width=0.44\textwidth]{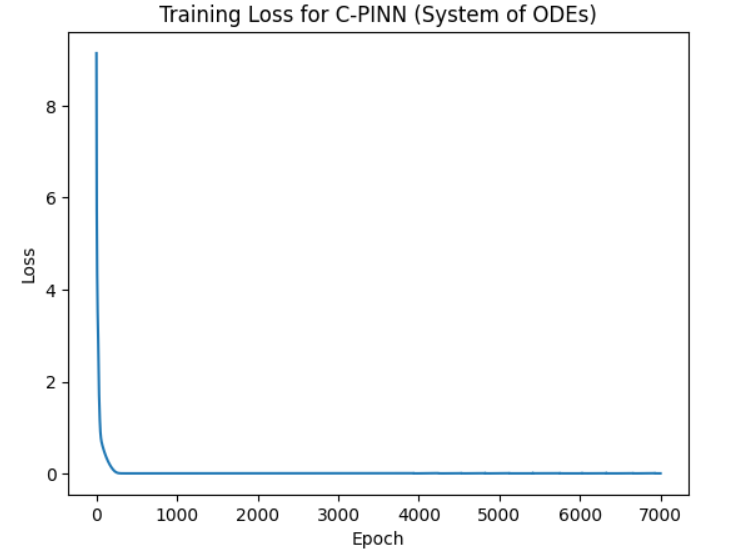}
	\includegraphics [width=0.46\textwidth]{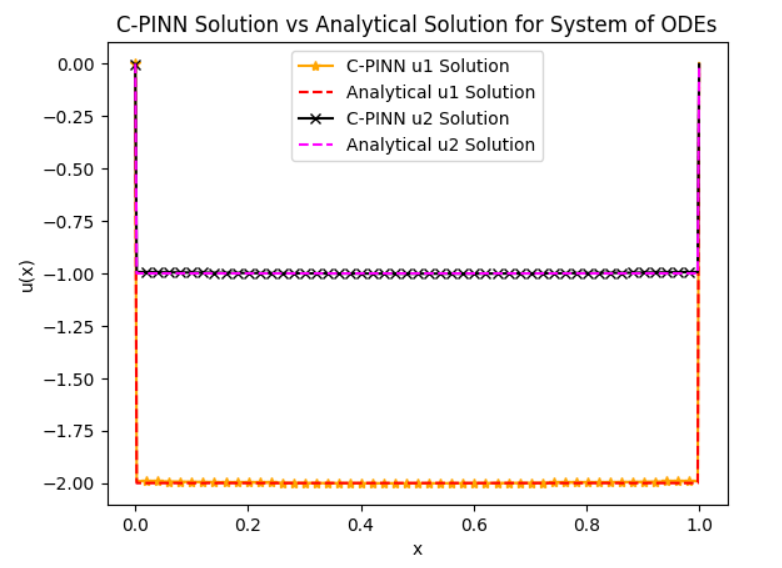}
	\caption{Loss and solution plot for C-PINN for different epochs.}
	\label{figure7}
\end{figure}

\begin{figure}[h]
	\centering
	\includegraphics [width=0.45\textwidth]{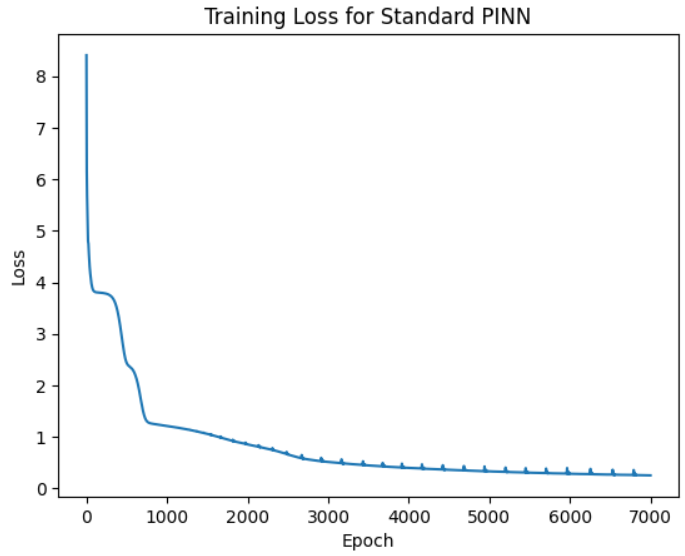}
	\includegraphics [width=0.45\textwidth]{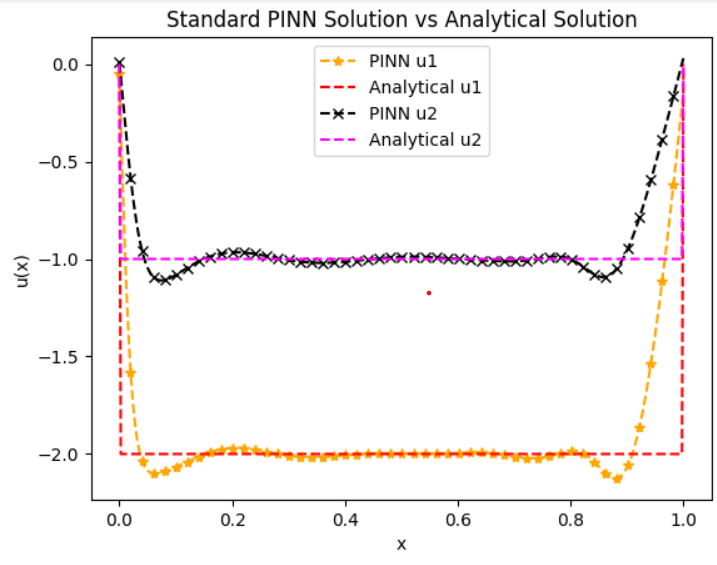}
	\caption{Loss and solution plot for PINN for different epochs.}
	\label{figure8}
\end{figure}

\begin{table}[ht]
	\footnotesize
	\centering
	\renewcommand{\arraystretch}{1}
	\begin{tabular}{|c|c|c|c|c|c|c|c|}
		\hline
		\textbf{Epoch} & 0 & 500 & 1000 & 1500 & 2000 & 2500 & 3000 \\
		\hline
		\textbf{C-PINN Loss} 
		& 9.144341 & 0.003289 & 0.003288 & 0.003301 & 0.003288 & 0.003559 & 0.003288 \\
		\textbf{PINN Loss} 
		& 8.406951 & 2.433742 & 1.212732 & 1.058357 & 0.854562 & 0.661478 & 0.514417 \\
		\hline
		\textbf{Epoch} & 3500 & 4000 & 4500 & 5000 & 5500 & 6000 & 6500 \\
		\hline
		\textbf{C-PINN Loss} 
		& 0.003288 & 0.003293 & 0.003287 & 0.003287 & 0.003288 & 0.003287 & 0.003287 \\
		\textbf{PINN Loss} 
		& 0.447261 & 0.398661 & 0.362370 & 0.331953 & 0.306092 & 0.285853 & 0.268869 \\
		\hline
	\end{tabular}
	\caption{Comparison of C-PINN loss and PINN loss at different epochs.}
	\label{table5}
\end{table}

Table \ref{table5} shows that initially,both PINN and C-PINN starts with comparable loss. However, C-PINN loss rapidly decreases to 0.003289 by epoch 500 and further stabilizes around 0.003287 at epoch 6500. Moreover, consistently achieves a lower loss compared to PINN across all epochs, indicating better convergence and accuracy. This comparison highlights the efficiency and stability of C-PINN in solving the given problem.
\newpage 
\subsection{ SPPs in two Dimensions (2D)} In this subsection, we present examples of two-dimensional singular perturbation problems, which are solved using the C-PINN method to obtain their solutions.
	
	\begin{algorithm}
		\caption{C-PINN Algorithm for 2D Singularly Perturbed Problems}
		\begin{algorithmic}[1]
			\State \textbf{Initialize:}
			\State \hspace{1em} Generate 2D collocation points $(x, y) \in [0,1]^2$ using LHS.
			\State \hspace{1em} Sample boundary points on $\partial\Omega$
			\State \hspace{1em} Set perturbation parameter $\varepsilon$, learning rate $\eta$, and max iterations $M$
					\State \hspace{1em} Initialize neural network weights: $\theta_1^o$, $\{\theta_1^k\}_{k=1}^K$, $\theta_2^o$, $\{\theta_2^k\}_{k=1}^K$ for $u_1$ and $u_2$
					\State Find the boundary layer location  to find $p(x)$ and $p(y)$
						\State Figure out the   boundary layer thickness to find  $\delta(\epsilon)_x$ and $\delta(\epsilon)_y$ 
			\State Define \textbf{OuterNN} for smooth interior solution
			\State Define \textbf{InnerNN}${x_i}$ and \textbf{InnerNN}${y_i}$ for layers at $x=i$ and $y=i$
			\State Construct composite network $u^\theta(x,y)$:
			\Statex \hspace{1.5em} $u^\theta(x,y) \gets \text{OuterNN}(x,y) + \text{InnerNN}_{x_i}(x,y)\, e^{-p(x)/\delta(\epsilon)_x} + \text{InnerNN}_{y_i}(x,y)\, e^{-p(y)/\delta(\epsilon)_y}$
			\State Define residual operator:
			\Statex \hspace{1.5em} $R(x,y) \gets \mathcal{D}^{\varepsilon}[u^\theta](x,y) - f(x,y)$
			\Statex \hspace{1.5em} (e.g., $\mathcal{D}^{\varepsilon}[u] = -\varepsilon \Delta u - u_x - u_y + u$)
			
			\For{$m = 1$ to $M$}
			\State Sample mini-batch of collocation and boundary points
			\State Compute $u^\theta(x,y)$ and its derivatives using auto-diff
			\State Evaluate residual $R(x,y)$ at collocation points
			\State Compute loss:
			\Statex \hspace{1.5em} $\mathcal{L}_{\text{res}} \gets \mathbb{E}_{(x,y)}[R(x,y)^2]$
			\Statex \hspace{1.5em} $\mathcal{L}_{\text{bc}} \gets \mathbb{E}_{(x,y) \in \partial\Omega} [u^\theta(x,y)^2]$
			\Statex \hspace{1.5em} $\mathcal{L} \gets \mathcal{L}_{\text{res}} + \mathcal{L}_{\text{bc}}$
			\State Update network parameters: $\theta \gets \theta - \eta \nabla_\theta \mathcal{L}$
			\EndFor
			
			\State \textbf{Return:} Trained composite model $u^\theta(x,y)$ as the 2D solution
			
			\State Use current weights as initialization for next step and continue the process for $M$ iterations.

		\end{algorithmic}
	\end{algorithm}

	\begin{example}\label{example2}
	Consider a singular perturbed convection-diffusion-reaction problem given below as:
	\begin{align}
		- \varepsilon \Delta u - (2 - x) u_x - u_y + \frac{3}{2} u = f \quad \text{in } \Omega = (0,1) \times (0,1)
	\end{align}
	with the boundary conditions:
	$$u = 0 \quad \text{on } \partial \Omega$$
	where 0$ < \epsilon \ll$ 1 is a small perturbation parameter.
	The source term f(x,y) is chosen such that the problem has the following analytical solution:
	$$u(x, y) = \left[\cos\left( \frac{\pi x }{2} \right) - \left(  \frac{e^{-x/\varepsilon} - e^{-1/\varepsilon}}{1 - e^{-1/\varepsilon}} \right) \right]
	\cdot \left( \frac{1 - e^{-y/\varepsilon}}{1 - e^{-1/\varepsilon}} \right) 
	$$
	\end{example}
	
In this 2D problem, the theoretical boundary layer is present at $x = 0$ and $y = 0$, with a thickness of order $\mathcal{O}(\epsilon).$ For our setup, we take $\epsilon = 10^{-5}$.
Our model struggled near $y=1,$ despite no theoretical boundary layer there. To address this, we added an auxiliary inner network around $y=1$, not to represent a physical layer, but to enhance flexibility and improve accuracy. Through empirical observations during training, this adjustment reflects the composite PINN strategy of using localized sub networks to refine solution accuracy and representation. The modification reduced global relative error and led to more consistent predictions across the domain. Now, let's go through the step-by-step implementation of the C-PINN method and the training process of the model.
	
First, we define two types of sub-neural networks: outer NNs and inner NNs. The outer networks capture the smooth variations, while the inner networks handle the sharp gradients of the coupled system. Each network consists of three fully connected layers with 100 and 150 neurons per layer respectively.

Next, we define a safe exponential function that clamps input values between -20 and 20 before applying the exponential operation, preventing overflow and underflow issue.
Furthermore, we create model instances, namely outer-nn-1, inner-nn-x0, inner-nn-y0 and inner-nn-u1 and utilize the following composite NN to represent a uniform solution:
{\scriptsize 	\begin{align}\nonumber
			\displaystyle & u(xy,\theta_1,\theta_2,\theta_3,\theta_4) =
			\text{outer\_nn}(xy,\theta_1) + \text{inner\_nn}_{x0}(xy,\theta_2) \cdot \text{safe\_exp}\left(\frac{-x}{\epsilon}\right)  \\ \nonumber
			& +  \text{inner\_nn}_{y0},(xy,\theta_3) \cdot \text{safe\_exp}\left(\frac{-y}{\epsilon}\right) 
			+\text{inner\_nn}_{y1}(xy,\theta_4) \cdot \text{safe\_exp}\left(\frac{-(1-y)}{\epsilon}\right). 
\end{align}}
We also define a function to calculate f(x,y) using python library SymPy and evaluate f(x,y) at the collocation points. 
Next, we define the loss function for C-PINN using PyTorch's built-in automatic differentiation. The loss formulation also includes the residual, ensuring the network learns to satisfy the governing equations. 
Lastly, we take 600 collocation points and  train our model for 7000 epochs with learning rate = $5*10^{-4}$ and visualize the results by plotting the loss curve and the solution graph. The solution graph compares the analytical solution of the coupled system with the C-PINN solution.

\begin{figure}[h]
	\centering
	\includegraphics [width=0.45\textwidth]{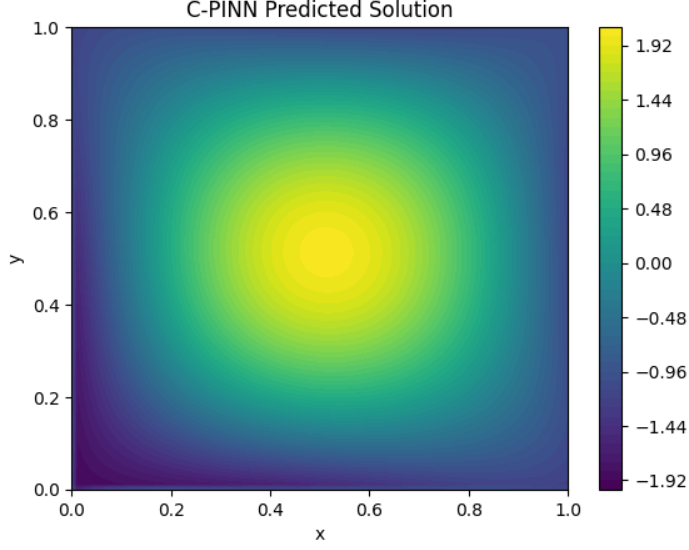}
	\includegraphics [width=0.45\textwidth]{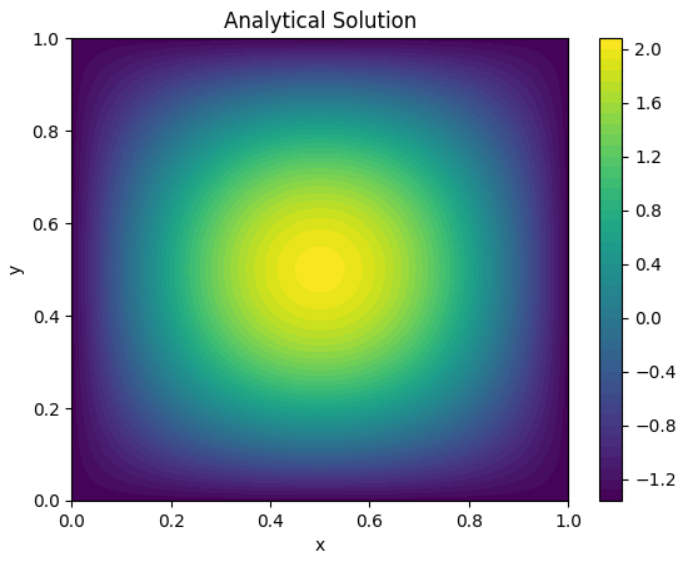}
	\caption{ C-PINN and analytical solution at different epochs.}
	\label{figure9}
\end{figure}

\begin{figure}[h]
	\centering
	\includegraphics [width=0.45\textwidth]{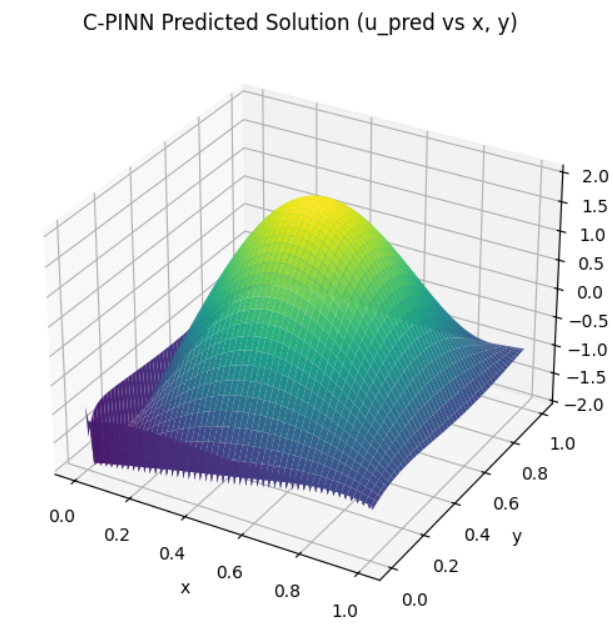}
	\includegraphics [width=0.45\textwidth]{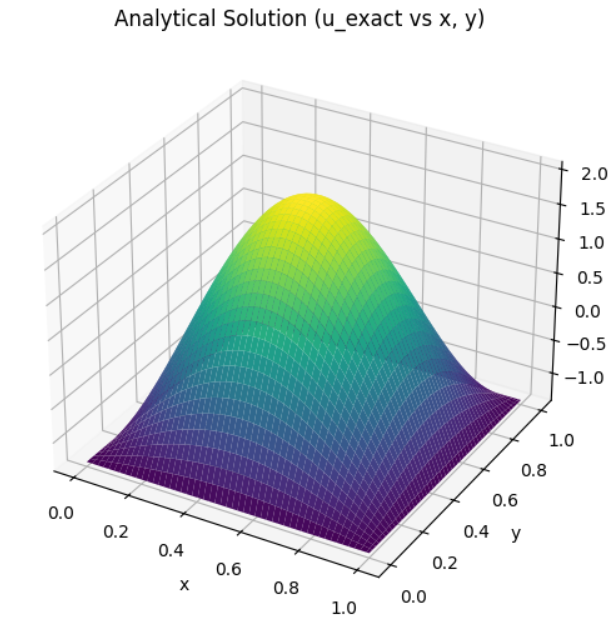}
	\caption{ C-PINN and analytical solution at different epochs.}
	\label{figure10}
\end{figure}	

\begin{figure}[h]
	\centering
	\includegraphics [width=0.45\textwidth]{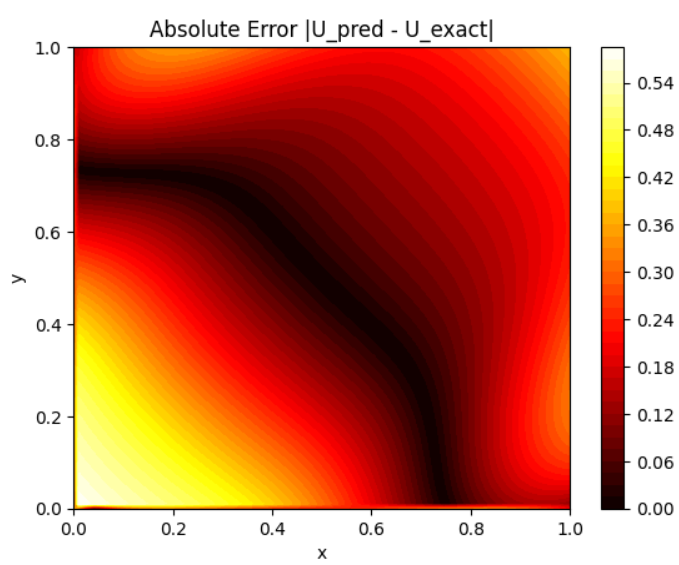}
	\includegraphics [width=0.45\textwidth]{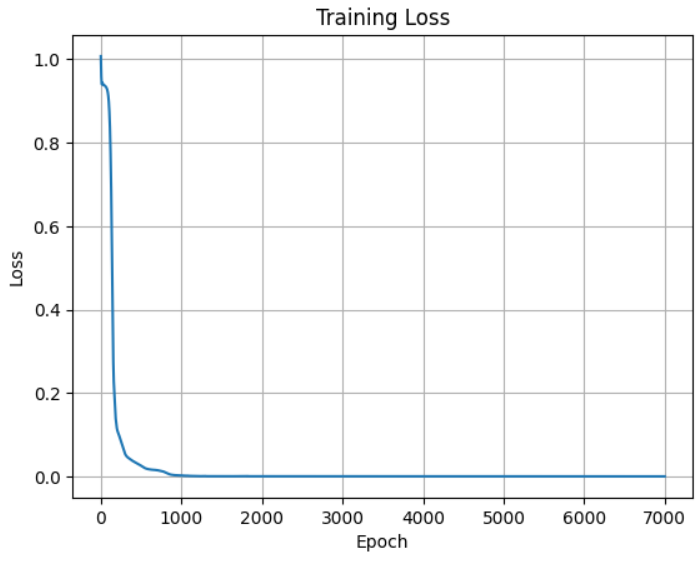}
	\caption{Loss plot for C-PINN for different epochs and absolute error plot.}
	\label{figure11}
\end{figure}
	
	Table \ref{table5} shows that initially, C-PINN starts with a loss of 5.785950. However, C-PINN loss rapidly decreases to 0.007740 by epoch 500 and further stabilizes around 0.000279 at epoch 6500. Moreover, consistently achieves a low loss across all epochs, indicating a good  convergence and accuracy.\\
	
	\begin{table}
		\centering
		\renewcommand{\arraystretch}{1.2}
		\setlength{\tabcolsep}{5pt}
		\footnotesize
		\begin{tabular}{|c|c | c |c |c |c |c |c|}
			\hline
			\textbf{Epoch}  & 0 & 500 & 1000 & 1500 & 2000 & 2500 & 3000  \\ 
			\hline
			\textbf{C-PINN Loss} & 5.785950 & 0.007740 & 0.004609 & 0.002732 & 0.001673 & 0.001294 & 0.000967  \\ 
			\hline
			\textbf{Epoch} & 3500 & 4000 & 4500 & 5000 & 5500 & 6000 & 6500 \\ 
			\hline
			\textbf{C-PINN Loss} & 0.000714 & 0.000549 & 0.000483 & 0.000418 & 0.000444 & 0.000316 & 0.000279  \\ 
			\hline
		\end{tabular}
		\caption{\( L_2 \) loss for C-PINN at different epochs.}
		\label{table5}
	\end{table}

\begin{example}\label{example3}

Consider a singular perturbed convection-diffusion-reaction problem given below as:
\begin{align}
	- \varepsilon \Delta u - u_x - u_y + u = f(x,y) \quad \text{in } \Omega = (0,1) \times (0,1),
\end{align}
with the boundary conditions:
$u = 0 \quad \text{on } \Gamma = \partial \Omega.$
where 0$ < \epsilon \ll$ 1 is a small perturbation parameter.
The source term f(x,y) is chosen such that the problem has the following analytical solution:
$ u(x, y) = \sin \pi x \cdot \sin \pi y \cdot \left(1 - e^{-x/\varepsilon}\right) \cdot \left(1 - e^{-y/\varepsilon}\right)$
\end{example}

In this 2D problem, the theoretical boundary layer is present at x = 0 and y = 0, with a thickness of order $\mathcal{O}(\epsilon).$ For our setup, we take $\epsilon = 10^{-5}$.
First, we define two types of sub-neural networks: outer NNs and inner NNs. The outer networks capture the smooth variations, while the inner networks handle the sharp gradients of the coupled system. Each network consists of three fully connected layers with 100 and 150 neurons per layer respectively.
Moreover, we create model instances, namely outer-nn-1, inner-nn-x0, inner-nn-y0 and inner-nn-y1 and utilize the following composite neural network to represent a uniform solution:
{\small 	\begin{align}\nonumber
			\displaystyle u(xy,\theta_1,\theta_2,\theta_3) & =
			\text{outer\_nn}(xy,\theta_1) + \text{inner\_nn}_{x0}(xy,\theta_2) \cdot \text{safe\_exp}\left(\frac{-x}{\epsilon}\right)  \\ \nonumber
			& +  \text{inner\_nn}_{y0},(xy,\theta_3) \cdot \text{safe\_exp}\left(\frac{-y}{\epsilon}\right). 
	\end{align}}

	\begin{figure}[h]
	\centering
	\includegraphics [width=0.45\textwidth]{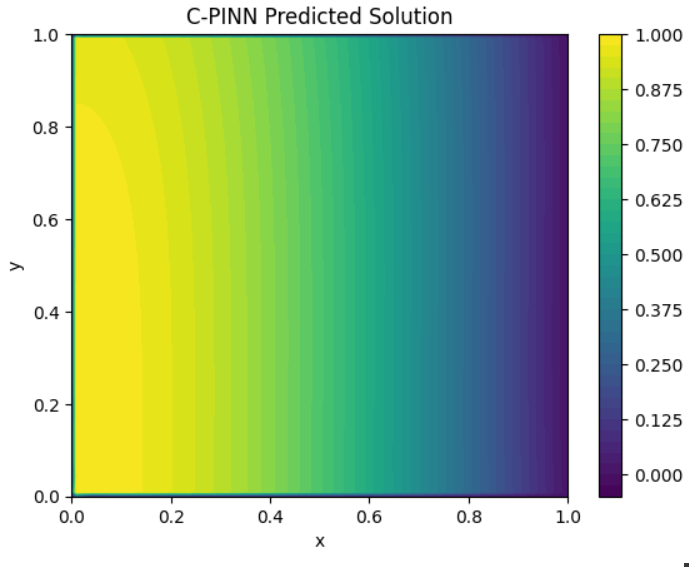}
	\includegraphics [width=0.45\textwidth]{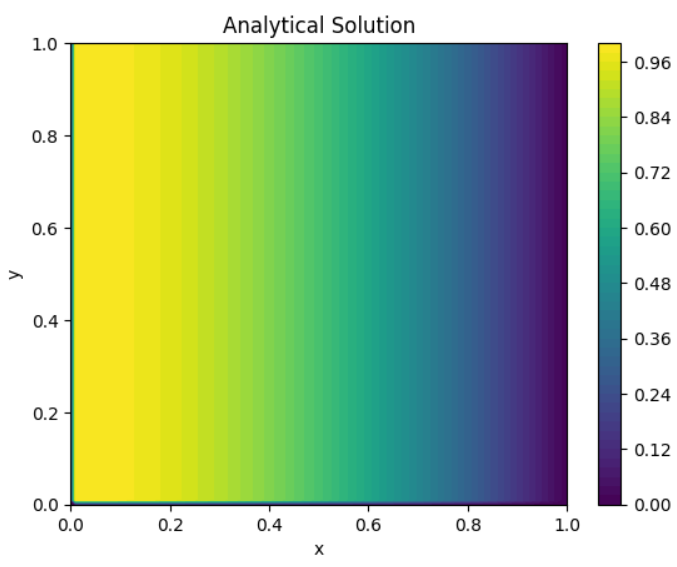}
	\caption{ C-PINN and Analytical Solution at different Epochs.}
	\label{figure12}
\end{figure}

\begin{figure}[h]
	\centering
	\includegraphics [width=0.45\textwidth]{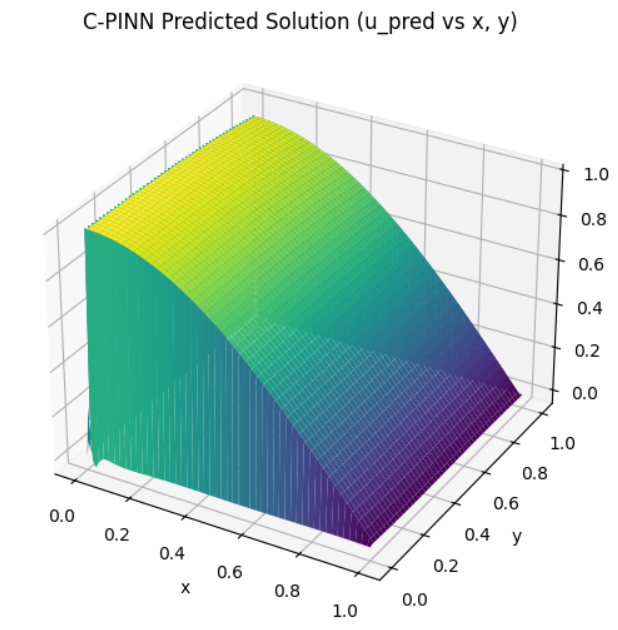}
	\includegraphics [width=0.45\textwidth]{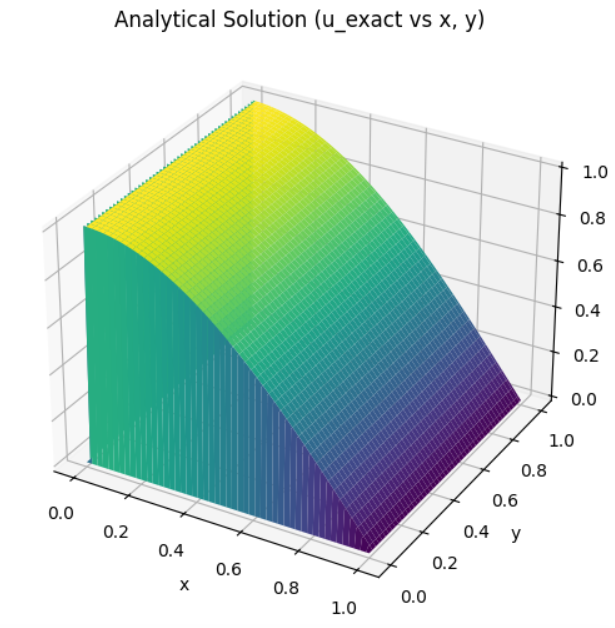}
	\caption{ C-PINN and Analytical Solution at different Epochs.}
	\label{figure15}
\end{figure}

\begin{figure}[h]
	\centering
	\includegraphics [width=0.45\textwidth]{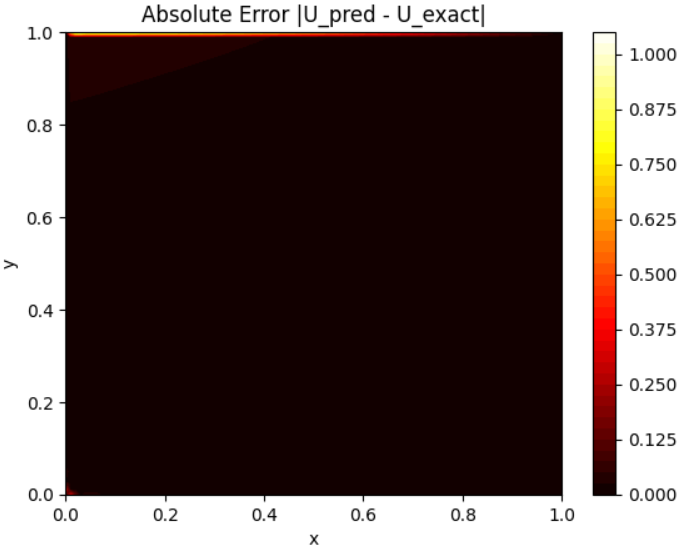}
	\includegraphics [width=0.45\textwidth]{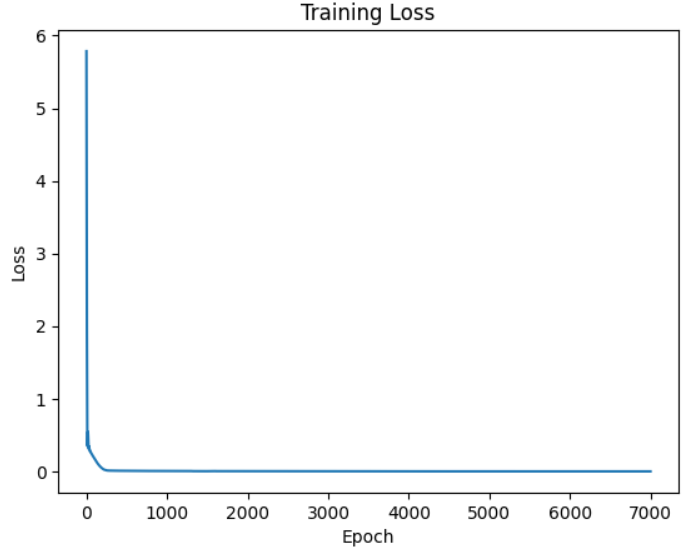}
	\caption{Loss plot for C-PINN for different epochs and absolute error plot.}
	\label{figure13}
\end{figure}

Table \ref{table6}  shows that initially, C-PINN starts with a loss of 1.007349. However, C-PINN loss rapidly decreases to 0.025603 by epoch 500 and further stabilizes around 0.000108 at epoch 6500. Moreover, consistently achieves a low loss across all epochs, indicating a good convergence and accuracy.

	\begin{table}[ht]
		\centering
		\renewcommand{\arraystretch}{1.4}
		\setlength{\tabcolsep}{5pt}
		\footnotesize
		\begin{tabular}{|c|c | c |c |c |c |c |c|}
			\hline
			\textbf{Epoch}  & 0 & 500 & 1000 & 1500 & 2000 & 2500 & 3000  \\ 
			\hline
			\textbf{C-PINN Loss} & 1.007349 & 0.025603 & 0.002031 & 0.000392 & 0.000333 & 0.000261 & 0.000217  \\ 
			\hline
			\textbf{Epoch} & 3500 & 4000 & 4500 & 5000 & 5500 & 6000 & 6500 \\ 
			\hline
			\textbf{C-PINN Loss} & 0.000180 & 0.000146 & 0.000132 & 0.000120 & 0.000113 & 0.000109 & 0.000108  \\ 
			\hline
		\end{tabular}
		\caption{\( L_2 \) loss for C-PINN at different epochs.}
		\label{table6}
	\end{table}

\newpage
	
	\section{Conclusion and Future Scope}\label{con}
	
	In this article, we studied the application of the C-PINN method to solve a wide range of SPPs, including convection-diffusion, reaction-diffusion and coupled systems of these equations in both one and two dimensions. To the best of our knowledge, this is the first work that applies C-PINNs to such a wide range of SPPs.
	Through several numerical experiments, we showed that the C-PINN method performs very well in capturing sharp boundary layers that typically arise in SPPs. Comparative analysis with standard PINNs confirms that C-PINN offers more accurate solutions. These results establish C-PINN as a more powerful method for solving SPPs.
	In the future, we plan to improve the method by developing adaptive techniques that can automatically find boundary layers and suggest the best neural network structure. We also aim to extend the C-PINN method to solve three-dimensional SPPs.

	\section*{Funding}
	The authors received no external funding for this work.

	\section*{Conflict of Interest}
	The authors declare that they have no conflicts of interest.
	
	\section*{Data Availability}
	Not applicable.

\end{document}